\documentclass[centertags,12pt]{amsart}
\usepackage{latexsym}
\usepackage{amssymb}

\numberwithin{equation}{section}
\makeatletter
\renewcommand{\subsection}{\@startsection
{subsection}{2}{0mm}{\baselineskip}{-0.25cm}
{\normalfont\normalsize\bf}}
\makeatother

\newtheorem{theorem}{Theorem}[section]
\newtheorem{proposition}[theorem]{Proposition}
\newtheorem{corollary}[theorem]{Corollary}
\newtheorem{lemma}[theorem]{Lemma}
{\theoremstyle{remark}
\newtheorem{remark}[theorem]{Remark}
\newtheorem{claim}{Claim}}
\theoremstyle{definition}

\def\P{\mathbf P}
\def\N{\mathbf N}
\def\Z{\mathbf Z}

\def\cF{\mathcal F}
\def\cH{\mathcal H}

\def\cX{\mathcal X}

\def\fp{\mathbf F_p}
\def\fq{\mathbf F_{q^2}}

\def\div{{\rm div}}

\def\deg{{\rm deg}}

\begin{document}
\author[Aguglia]{Angela Aguglia}
\author[Korchm\'aros]{G\'abor Korchm\'aros}
\author[Torres]{Fernando Torres}
\thanks{This research was carried out within the project
``Progetto e Realizzazione di un Criptosistema per
Telecomunicazioni'', POP FESR 1994/99 -- II Triennio. The third
author was partially supported by Cnpq-Brazil, Proc. 300681/97-6}
\title[Plane maximal curves]{Remarks on plane maximal curves}
\address{Dipartimento di Matematica, Universit\`a di Basilicata\\
Via N. Sauro 85, 85100 Potenza, Italy}
\email{aguglia@matna2.dma.unina.it}
\email{korchmaros@unibas.it}
\address{IMECC-UNICAMP, Cx. P. 6065, Campinas, 13083-970-SP, Brazil}
\email{ftorres@ime.unicamp.br}

     \begin{abstract} Some new results on plane $\fq$-maximal
curves are stated and proved. By \cite{r-sti}, the degree $d$ of a
plane $\fq$-maximal curve is less than or equal to $q+1$ and
equality holds if and only if the curve is $\fq$-isomorphic to the
Hermitian curve. We show that $d\leq q+1$ can be improved to $d\leq
(q+2)/2$ apart from the case $d=q+1$ or $q\leq 5$. This upper bound
turns out to be sharp for $q$ odd. In \cite{carbonne-henocq} it was
pointed out that some Hurwitz curves are plane $\fq$-maximal
curves. Here we prove that (\ref{eq1.2}) is the necessary and
sufficient condition for a Hurwitz curve to be $\fq$-maximal. We
also show that this criterium holds true for the $\fq$-maximality
of a wider family of curves.
   \end{abstract}

\maketitle

   \section{Introduction}\label{s1}

A $\fq$-maximal curve of genus $g$ is a projective, geometrically
irreducible, non-singular, algebraic curve defined over a
finite field $\fq$ of order $q^2$ such that the number of its
$\fq$-rational points attains the Hasse-Weil upper bound
   $$
 1+q^2+2qg\, .
   $$
Maximal curves, especially those having large genus with respect to
$q$, are known to be very useful in Coding theory \cite{goppa}.
Also, there are various ways of employing them in Cryptography, and
it is expected that this interesting connection will be be explored
more fully, see \cite[Chapter 8]{shpar}. Other motivation for the
study of maximal curves comes from Correlations of Shift Register
Sequences \cite{L-N}, Exponentials Sums over Finite Fields
\cite{moreno}, and Finite Geometry \cite{hirschfeld}. Recent papers
on maximal curves which also contain background and expository
accounts are \cite{r-sti}, \cite{sti-x}, \cite{ft}, \cite{fgt},
\cite{geer-vl}, \cite{ft1}, \cite{ckt1}, \cite{gt}, \cite{chkt},
\cite{at}, \cite{ckt2}, and \cite{kt}.

A relevant result on $\fq$-maximal curves $\cX$ with genus $g$ states that
either $g=q(q-1)/2$ and $\cX$ is $\fq$-isomorphic to the Hermitian curve
$\cH$ of equation
   \begin{equation}\label{hermitian}
X^{q+1}+Y^{q+1}+Z^{q+1}=0\, ,
  \end{equation}
or $g\leq(q-1)^2/4$; see \cite{ihara}, \cite{sti-x}, and \cite{ft}.
One expects that the bound $(q-1)^2/4$ can be substantially lowered
apart from a certain number of exceptional values of $g$. Finding
such values is one of the problems of current interest in the study
of maximal curves; see \cite[Section 3]{fgt}, \cite[Proposition
2.5]{ft1}, \cite[Section 3]{ckt1}, and \cite{at}.

In this paper we investigate plane maximal curves. In Section
\ref{s2} we prove the non-existence of a plane $\fq$-maximal curve
whose genus belongs to the interval $(q(q-2)/8,q(q-2)/4]$, for $q$
even, and $((q-1)(q-3)/8,(q-1)^2/4]$ for $q$ odd; see Corollary
\ref{cor2.1}. The curves studied in Section \ref{s3} show that
these bounds are sharp in some cases. In contrast, a few examples
of (non planar) $\fq$-maximal curves with genera in these intervals
are known to exist; see \cite[Section 3]{fgt}, \cite[pp.
74--75]{ckt1}, \cite{at}, \cite{g-sti-x}, and \cite[Theorem
2.1]{ckt2}.

In the course of our investigation we point out that the Hermitian
curve $\cH$ is the unique $\fq$-maximal curve (up to
$\fq$-isomorphism) which is $\fq$-Frobenius non-classical with
respect to the linear series $\Sigma_1$ cut out by lines; see
Proposition \ref{prop2.1}. Also, the order of contact $\epsilon_2$
of a non-classical (with respect to $\Sigma_1$) $\fq$-maximal curve
with the tangent at a general point satisfies $\epsilon_2^2\le
q/p$, where $p:={\rm char(\fq)}$; see Corollary \ref{cor2.2}. In
particular, plane $\fq$-maximal curves with $q=p$ and $q=p^2$ are
classical with respect to $\Sigma_1$.

According to \cite[Prop. 6]{lachaud}, every curve which is $\fq$-covered
by the Hermitian curve is $\fq$-maximal. An open problem of considerable
interest is to decide whether the converse of this statement also holds.
In Section \ref{s3} we solve this problem for the family of the so-called
Hurwitz curves. Recall that a Hurwitz curve of degree $n+1$ is defined as
a non--singular plane curve of equation
  \begin{equation}
  \label{hurwitz}
X^nY+Y^nZ+Z^nX=0\, ,
  \end{equation}
where $p={\rm char(\fq)}$ does not divide $n^2-n+1$. Theorem \ref{thm3.1}
together with Corollary \ref{cor3.1} states indeed that the Hurwitz curve
is $\fq$-covered by the Hermitian curve if and only if
        \begin{equation}\label{eq1.2}
q+1\equiv0\pmod{(n^2-n+1)}\,.
        \end{equation}

It should be noted on the other hand that for certain $n$ and $p$,
the Hurwitz curve is not $\fq$-maximal for any power $q$ of $p$;
this occurs, for instance, for $n=3$ and $p\equiv 1 \pmod{7}$. One
can then ask for conditions in terms of $n$ and $p$ which assure
that the Hurwitz curve is $\fq$-maximal for some power $q$ of $p$.
Our results in this direction are given in Remarks \ref{rem3.0} and
\ref{rem3.01}, and Corollaries \ref{cor3.3}, \ref{cor3.4}. They
generalise some previous results obtained in \cite[Lemmes 3.3,
3.6]{carbonne-henocq}. Another feature of the Hurwitz curve is that
it is non-classical provided that $p^e$ divides $n$ with $p^e\ge
3$; see Remark \ref{rem3.1}. So if both (\ref{eq1.2}) and $p^e|n$
hold then the Hurwitz curve turns out to be a non-classical plane
$\fq$-maximal curve. As far as we know, these Hurwitz curves
together with the Hermitian curves and the Fermat curves of degree
$n^2-n+1$ (see Corollary \ref{cor3.1}), are the only known examples
of non-classical plane $\fq$-maximal curves. As mentioned before,
these curves show the sharpness of some of the results obtained in
Section \ref{s2}. 

Hurwitz curves as well as their generalizations have been investigated for
several reasons by many authors; see \cite[Section 1]{bennama-carbonne}
and \cite{pellikaan}. This gives a motivation to the final Section
\ref{s4} where we show that the main results of Section \ref{s3} extend to
(the non--singular model of) the curve with equation
  $$
X^nY^\ell+Y^nZ^\ell+Z^nX^\ell=0\, ,
  $$
where $n\ge \ell\ge 2$ and $p={\rm char}(\fq)$ does not divide
$Q(n,\ell):= n^2-n\ell+\ell^2$.

Our investigation uses some concepts, such as non-classicity, from
St\"ohr-Voloch's paper \cite {sv} where an alternative proof to the
Hasse-Weil bound was given among other things. We also refer to
that paper for terminology and background results on orders and
Frobenius orders of linear series on curves.

    \section{The degree of a plane maximal curve}\label{s2}

Let $\cX$ be a plane $\fq$-maximal curve of degree $d\geq 2$. Since
the genus of $\cX$ is equal to $(d-1)(d-2)/2$, the upper bound for
$g$ quoted in Sec. \ref{s1} can be rephrased in terms of $d$:
  \begin{equation}\label{eq1.1}
d\leq d_1(q):=\frac{3+\sqrt{2(q-3)(q+1)+9}}{2}   \qquad\text{or}\qquad
d=q+1\, .
  \end{equation}
The main result in this section is the improvement of (\ref{eq1.1})
given in Theorem \ref{thm2.1}: Apart from small $q$'s, either
$d=q+1$, or $d=\lfloor (q+2)/2\rfloor$, or $d$ is upper bounded by
a certain function $d_5(q)$ such that $d_5(q)/q\approx 2/5$.
Our first step consists in lowering $d_1(q)$ to $d_2(q)$ with $d_2(q)/q
\approx 1/2$.

Let $\Sigma_1$ be the linear series cut out by lines of
$\P^2(\bar\fq)$ on $\cX$. For $P\in \cX$, let $j_0(P)=0<j_1(P)=1<j_2(P)$
be the $(\Sigma_1,P)$-orders, and $\epsilon_0=0<\epsilon_1=1<\epsilon_2$
(resp. $\nu_0=0<\nu_1$)  the orders (resp. $\fq$-Frobenius orders) of
$\Sigma_1$. We let $p$ be the characteristic of $\fq$.
   \begin{lemma}\label{lemma2.1}
  \begin{enumerate}
\item[\rm(1)] $\nu_1\in \{1,\epsilon_2\};$
\item[\rm(2)] $\epsilon_2\leq q;$
\item[\rm(3)] $\epsilon_2$ is a power of $p$ whenever $\epsilon_2>2.$
  \end{enumerate}
  \end{lemma}
     \begin{proof} For (1), see \cite[Prop.
2.1]{sv}. For (2), suppose that $\epsilon_2>q$, then $\epsilon_2=q+1$
as $\epsilon_2\leq d$ and $d\le q+1$ by (\ref{eq1.1}). Then, by the
$p$-adic criterion \cite[Cor. 1.9]{sv}, $q$ would be a $\Sigma_1$-order,
a contradiction. For (3), see \cite[Prop. 2]{garcia-voloch}.
     \end{proof}
The following result is a complement to \cite[Prop. 3.7]{pardini},
\cite[Thm. 6.1]{homma}, and \cite[Prop. 6]{hefez-voloch}.
      \begin{proposition}\label{prop2.1} For a plane
$\fq$-maximal curve $\cX$ of degree $d\ge 3$, the following conditions are
equivalent:
  \begin{enumerate}
\item[\rm(1)] $d=q+1;$
\item[\rm(2)] $\cX$ is $\fq$-isomorphic to the Hermitian
of equation (\ref{hermitian})$;$
\item[\rm(3)] $\epsilon_2=q;$
\item[\rm(4)] $\nu_1=q;$
\item[\rm(5)] $j_2(P)=q+1$ for each $P\in \cX(\fq);$
\item[\rm(6)] $\nu_1>1$; i.e, $\Sigma_1$ is $\fq$-Frobenius
non-classical$.$
  \end{enumerate}
  \end{proposition}
  \begin{proof} $\text{(1)$\Rightarrow$(2)}:$ Since the genus
of a non--singular plane curve of degree $d$ is $q(q-1)/2$, part
(2) follows from \cite{r-sti}.

$\text{(2)$\Rightarrow$(3)}:$ This is well known property of the
Hermitian curve; see e.g. \cite[p. 105]{ft} or \cite{garcia-viana}.

$\text{(3)$\Rightarrow$(4)}:$ If $q=2$, then from $d\ge
\epsilon_2=q$ and (\ref{eq1.1}), either $d=2$ or $d=3$. By hypothesis,
$d=3$ can only occur, and so, by parts (1) and (2), $\cX$ is
$\mathbf F_4$-isomorphic to the Hermitian curve $X^3+Y^3+Z^3=0$.
Then $\nu_1=\epsilon_2=2$; see loc. cit.

Let $q\ge 3$. By Lemma \ref{lemma2.1}(1), $\nu_1\in\{1,q\}$.
Suppose that $\nu_1=1$ and let $S_1$ be the $\fq$-Frobenius divisor
associated to $\Sigma_1$. Then \cite[Thm. 2.13]{sv}
   $$
\deg(S_1)=(2g-2)+(q^2+2)d \geq 2\#\cX(\fq)=2(q+1)^2+2q(2g-2)
   $$
so that $((2q-1)d-(q^2+2q+1))(d-2)\le 0$, and hence
   \begin{equation}\label{eq2.1}
d\le F(q):=(q^2+2q+1)/(2q-1)\, .
    \end{equation}
Thus, as $d\ge \epsilon_2=q$, we would have $q^2-3q-1\le 0$ and hence
$q\le 3$. If $q=3$, from (\ref{eq2.1}) we have that $d=3$; this
contradicts \cite[Cor. 2.2]{pardini} (cf. Remark \ref{rem2.2}(ii)).

$\text{(4)$\Rightarrow$(5)}:$ By \cite[Cor. 2.6]{sv}, $\nu_1 \leq
j_2(P)-1$ for any $P \in \cX(\fq)$. Then part (5) follows as $j_2(P) \leq
d$ and $d \leq q+1$ by (\ref{eq1.1}).

$\text{(5)$\Rightarrow$(6)}:$ Suppose that $\nu_1=1$. Then, by \cite[Prop.
2.4(a)]{sv}, $v_P(S_1)\ge q+1$ for any $P\in \cX(\fq)$. Therefore
$$
\deg(S_1)=(2g-2)+(q^2+2)d\ge (q+1)\#\cX(\fq)=(q+1)^3+(q+1)q(2g-2)\, ,
$$
a contradiction as $3\le d\le q+1$.

$\text{(6)$\Rightarrow$(1)}:$ From \cite[Thm. 1]{hefez-voloch} and the
$\fq$-maximality of $\cX$ we  have
    $$
\#\cX(\fq)=d(q^2-1)-(2g-2)=(1+q)^2+q(2g-2)\, .
    $$
Since $2g-2=d(d-3)$ and $d>1$, part (1) follows.
     \end{proof}
     \begin{corollary}\label{cor2.1} Let $d\ge 3$ be the degree of a
     plane $\fq$-maximal curve. Then either $d=q+1$ or
  \begin{equation*}
d\leq d_2(q):=\begin{cases}
\lfloor (q+2)/2\rfloor & \text{if $q\ge 4$ and $q\neq 3,5$,}\\
          3            & \text{if $q=3$,}\\
          4            & \text{if $q=5$.}
\end{cases}
   \end{equation*}
In particular, for $q\neq 3,5$, a $\fq$-maximal curve has no
non-singular plane model if its genus is assumed to belong to the
interval $(q(q-2)/8,q(q-2)/4]$, for $q$ even, and
$((q-1)(q-3)/8,(q-1)^2/4]$, for $q$ odd.
   \end{corollary}
   \begin{proof} The statement on the genus follows immediately from the
upper bound on $d$. By (\ref{eq1.1}) we have that $d\le q+1$. If
$d<q+1$, then $q\ge 3$ and from Proposition \ref{prop2.1}
$\Sigma_1$ is $\fq$-Frobenius classical. In particular,
(\ref{eq2.1}) holds true; i.e., we have $d\le F(q)$. It is easy to
see that $F(q)<(q+3)/2$ for $q>5$ and that $F(4)=25/7$. Moreover,
$F(3)=16/5$ and $F(5)=4$, and the result follows.
   \end{proof}
   \begin{remark}\label{rem2.1} Let $d$ be the degree of a plane
$\fq$-maximal curve of degree $d$ and assume that $3\le d\le
d_2(q)$.

(i) If $q$ is odd, then the $\fq$-maximal curve of equation
$$
X^{(q+1)/2}+Y^{(q+1)/2}+Z^{(q+1)/2}=0\, ,
$$
shows that the upper bound $d_2(q)=(q+1)/2$ in Corollary
\ref{cor2.1} is the best possible as far as $q\ne 3,5$. We notice
that this curve is the unique $\fq$-maximal plane curve (up to
$\fq$-isomorphism) of degree $(q+1)/2$ provided that $q\ge 11$; see
\cite{chkt}.

(ii) From results of Deuring, Tate and Watherhouse (see e.g. \cite[Thm.
4]{ughi}), there exist elliptic $\fq$-maximal curves for any $q$. In
particular, $d_2(q)=3$ is sharp for $q=3$.

(iii) From \cite[Sec. 4]{serre}, there exists a plane quartic $\mathbf
F_{25}$-maximal; so $d_2(q)=4$ is sharp for $q=5$.

(iv) By part (ii), $d_2(q)=3$ is sharp for $q=4$. For $q\ge 8$, $q$
even, no information is currently available to asses how good the
bound $d_2(q)=(q+2)/2$ is.
   \end{remark}
We go on to look for an upper bound for the degree $d$ of a
$\fq$-maximal curve satisfying the condition $d<\lfloor
(q+2)/2\rfloor$. Our approach is inspired on \cite[Sec. 3]{chkt}
where the $\fq$-Frobenius divisor $S_2$ associated with the linear
series $\Sigma_2$ cut out on $\cX$ by conics was employed to obtain
upper bounds for the number of $\fq$-rational points of plane
curves. In fact, if we use $\Sigma_2$ instead of $\Sigma_1$, we can
get better results for values of $d$ ranging in certain intervals
depending on $q$. This was pointed out at the first time in
\cite{garcia-voloch1}. \\ In order to compute the $\Sigma_2$-orders
of a plane $\fq$-maximal curve $\cX$, one needs to know whether
$\cX$ is classical or not with respect to $\Sigma_1$. This gives
the motivation to Proposition \ref{prop2.2}. The following remark
will be useful in the proof.
   \begin{remark}\label{rem2.2} (i) If a projective, geometrically
irreducible, non-singular, algebraic curve defined over a field of
characteristic $p>0$ admits a linear series $\Sigma$ of
degree $D$, then $\Sigma$ is classical provided that $p>D$; see \cite[Cor.
1.8]{sv}.

(ii) If a non-singular plane curve of degree $D$ defined over a field of
characteristic $p>$ is non-classical with respect to the linear series
cut out by lines, then $D\equiv 1\pmod{p}$; see \cite[Cor. 2.2]{pardini},
and \cite[Cor. 2.4]{homma1}.
   \end{remark}
   \begin{proposition}\label{prop2.2} Let $\cX$ be a plane $\fq$-maximal
curve of degree $d$ such that $3\le d\le d_2(q)$, where $d_2(q)$ is as in
Corollary \ref{cor2.1}. Then the linear series
$\Sigma_1$ on $\cX$ is classical provided that one of the following
conditions holds:
   \begin{enumerate}
\item[\rm(i)] $p\ge d$ or $d\not\equiv 1\pmod{p};$
\item[\rm(ii)] $q= 4,8,16,32;$
\item[\rm(iii)] $p\ge 3$ and either $q=p$ or $q=p^2;$
\item[\rm(iv)] $p=2$, $q\ge 64$, and either $d\le 4$, or $d\ge d_3(q):=
q/4-1$ for
$q=64, 128, 256$, or $d\ge d_3(q):=q/4$ for $q\ge 512;$
\item[\rm(v)] $p\ge 3$, $q=p^v$ with $v\ge 3$, and $d\ge d_3(q):=
q/p-p+2.$
   \end{enumerate}
   \end{proposition}
   \begin{proof} If (i) holds, then $\Sigma_1$ is classical by Remark
\ref{rem2.2}. For $q=p$, the hypothesis on $d$ yields $p\ge 3$ and
hence $d\le (p+1)/2<p$. Thus $\Sigma_1$ is classical by Remark
\ref{rem2.2}(i). Note that the following computations will provide
another proof of this fact. \\ For the rest of the proof we assume
$\Sigma_1$ to be non-classical, and we show that no one of the
conditions (i),\ldots,(v) holds. From Lemma \ref{lemma2.1}(3),
$\epsilon_2\ge M$, where $M=4$ for $p=2$, and $M=p$ for $p\ge 3$.
Also, $\nu_1=1$ by Proposition \ref{prop2.1}. Therefore, as
$j_2(P)\ge
\epsilon_2$ for each $P\in\cX$ \cite[p. 5]{sv}. From \cite[Prop.
2.4(a)]{sv} we have that $v_P(S_1)\ge M$ for each $P\in
\cX(\fq)$, where as before $S_1$ denotes the $\fq$-Frobenius divisor
associated to $\Sigma_1$. Thus,
   $$
\deg(S_1)=(2g-2)+(q^2+2)d\ge M\#\cX(\fq)=M(q+1)^2+Mq(2g-2)\, ,
   $$
or, equivalently,
   $$
(Mq-1)d^2-(q^2+3Mq-1)d+M(q+1)^2\le 0\, .
   $$
On the other hand, the discriminant of the above quadratic
polynomial in $d$ is
   $$
\Delta_M(q):=q^4-(4M^2-6M)q^3+(M^2+4M-2)q^2-(4M^2-2M)q+4M+1\, ,
   $$
and hence $\Delta_M(q)<0$ if and only if either $q=4,8,16,32$ and
$M=4$, or $q=p, p^2$ and $M=p\ge 3$. For these $q$'s, the above
inequality cannot actually hold, and hence $\Sigma_1$ must be
classical. Furthermore, if $\Delta_M(q)\ge 0$, then
$$
F'(M,q):=\frac{q^2+3Mq-1-\sqrt{\Delta_M(q)}}{2(Mq-1)}\le d\le
F(M,q):= \frac{q^2+3Mq-1+\sqrt{\Delta_M(q)}}{2(Mq-1)}\, .
$$
It is easy to check that $F'(4,q)>4$, $F(4,q)<q/4-1$ for $q=64,
128,256$, and that $F(4,q)<q/4$ for $q\ge 512$; hence if (iv)
holds, then $\Sigma_1$ must be classical. Let $p\ge 3$. If
$q/p-p+2\le d\le q/p$, then $\Sigma_1$ must be classical by (i). So
we can suppose that $d\ge q/p+1$. It turns out that $F(p,q)<q/p+1$
and hence the result follows when (v) is assumed to be true.
   \end{proof}
   \begin{remark}\label{rem2.3} For $q=p^3$, $p\ge 3$, the bound $d_3(q)$
in Proposition \ref{prop2.2} is sharp. Indeed, there exists a plane
$\mathbf F_{p^6}$-maximal curve of degree $p^2-p+1$ which is
non-classical for $\Sigma_1$; see Corollary \ref{cor3.1} and Remark
\ref{rem3.1}.
    \end{remark}
    \begin{corollary}\label{cor2.2} Let $\cX$ be a plane $\fq$-maximal curve of
degree $d$ as in Proposition \ref{prop2.2}. Assume that $\cX$ is
non-classical for
$\Sigma_1$ and let $\epsilon_2$ be the order of contact of $\cX$ with the
tangent at a general point. Then
   \begin{enumerate}
\item[\rm(1)] $q\ge 64$ if $p=2$, and $q\ge p^3$ for $p\ge 3;$
\item[\rm(2)] $\epsilon_2^2\le q/p.$
   \end{enumerate}
    \end{corollary}
    \begin{proof} Part (1) follows from Proposition
\ref{prop2.2}(ii)(iii). To prove (2), we first note that
$\epsilon_2<q$ (cf. Proposition \ref{prop2.1}), and that
$\epsilon_2$ is a power of $p$ (see Lemma \ref{lemma2.1}(3)). Now,
with the same notation as in the proof of the previous proposition,
we get $d\le F(M,q)$ with $M=\epsilon_2$. So $d\le q/\epsilon_2$.
Furthermore, $d\ge\epsilon_2$ and so $d\ge
\epsilon_2+1$ by Remark \ref{rem2.2}(ii). Hence $\epsilon_2+1\le
q/\epsilon_2$ and part (2) follows.
    \end{proof}
   \begin{remark}\label{rem2.3.1} The example in Remark \ref{rem2.3} shows
that Corollary \ref{cor2.2}(1) is sharp for $p\ge 3$.
   \end{remark}
Our next step is to show that every plane $\fq$-maximal curve which
is classical for $\Sigma_1$ contains an $\fq$-rational point
different from its inflexions.
   \begin{lemma}\label{lemma2.2} Let $\cX$ be a $\fq$-maximal curve
of degree $d\ge 3$ which
is classical with respect to $\Sigma_1$. Then there exists $P_0\in
\cX(\fq)$ whose $(\Sigma_1,P_0)$-orders are $0, 1,2.$
   \end{lemma}
   \begin{proof} Let $R_1$ be the ramification divisor associated to
$\Sigma_1$ and suppose that $j_2(P)\ge 3$ for each $P\in \cX(\fq)$. Then
from \cite[p. 12]{sv},
$$
\deg(R_1)=3(2g-2)+3d\ge \#\cX(\fq)=(q+1)^2+q(2g-2)
$$
which is a contradiction as $g\ge 1$ and $3\le d <q+1$.
   \end{proof}

It should be noticed that Lemma \ref{lemma2.2} improves a previous
result, see \cite[Cor. 3.2]{chkt}.

We are in a position to establish some useful properties of the
linear series $\Sigma_2$ cut out by conics of $\P^2(\bar\fq)$ on
plane $\fq$-maximal curve $\cX$ of degree $d\ge 3$. Since $\cX$ is
non-singular, $\Sigma_2=2\Sigma_1$. Taking into account $d\ge 3$,
we see that $\Sigma_2$ is a 5-dimensional linear series of degree
$2d$.
    \begin{lemma}\label{lemma2.3} Let $d$ be the degree of a plane
$\fq$-maximal curve $\cX$. Let $q=8$ or $q\ge 11$, and suppose that
   $$
d_4(q):=\frac{2q^2+15q-20+\sqrt{4q^4-40q^3+145q^2-300q+600}}{10(q-2)}<d
\le d_2(q)\, ,
   $$
where $d_2(q)$ is as in Corollary \ref{cor2.1}. Then the orders
(resp. $\fq$-Frobenius orders) of $\Sigma_2$ are
$0,1,2,3,4,\epsilon$ (resp. $0,1,2,3,\epsilon$) with
$5\le\epsilon\le q$. Furthermore, $p$ divides $\epsilon.$
   \end{lemma}

   \begin{proof} By some computations we obtain that $d_4(q)$ is bigger
than $d_3(q)$ in Proposition \ref{prop2.2}. So the curve $\cX$ is
classical for $\Sigma_1$. Let $P_0\in \cX(\fq)$ be as in Lemma
\ref{lemma2.2}. Then the $(\Sigma_2,P_0)$-orders are $0,1,2,3,4$
and $j_0$ with $5\le j_0\le 2d$ (cf. \cite[p. 464]{garcia-voloch}).
Therefore, the $\Sigma_2$-orders are $0,1,2,3,4$ and $\epsilon$
with $5\le \epsilon\le j_0$. Since $j_0\le 2d$, from Corollary
\ref{cor2.1}, $\epsilon\le q+2$, and hence $\epsilon\le q$ by the
$p$-adic criterion \cite[Cor. 1.9]{sv}. Also, the $\fq$-Frobenius
orders of $\Sigma_2$ are $0,1,2,3$ and $\nu$ with
$\nu\in\{4,\epsilon\}$; see \cite[Prop. 2.1, Cor. 2.6]{sv}. Suppose
that $\nu=4$ and keep up $S_2$ to denote the $\fq$-Frobenius
divisor associated to $\Sigma_2$. Then \cite[Thm. 2.13]{sv}
   $$
\deg(S_2)=10(2g-2)+(q^2+5)2d\ge 5\#\cX(\fq)=5(q+1)^2+5q(2g-2)
   $$
or equivalently
  $$
(5q-10)d^2-(2q^2+15q-20)d+5(q+1)^2\le 0\, .
  $$
The discriminant of this equation is $4q^2-40q^3+145q^2-300q+600$ and it
is positive for any $q$. Since $d_4(q)$ is the biggest root of the
quadratic polynomial in $d$ above, $d\le d_4(q)$, a contradiction.
Finally, $p$ divides $\epsilon$ by \cite[Cor. 3]{garcia-homma}.
   \end{proof}
Let $d_4(q)$ be as in Lemma \ref{lemma2.3} and for $q=p^v$, $v\ge
2$, let $d_4(p,q)$ denote the function
  $$
\frac{2q^2+3(5-\frac{1}{p})q-8+\sqrt{4q^4-8(5-\frac{1}{p})q^3+
(113-\frac{50}{p}+\frac{9}{p^2})q^2
-4(25-\frac{17}{p})q+184}}{2(5-\frac{1}{p})q-12}\, .
  $$
  \begin{theorem}\label{thm2.1} Let $d$ be the degree of a plane
$\fq$-maximal curve $\cX$. Suppose that $3\le d<q+1$ and that $q=8$ or
$q\ge 11$. Then
\begin{equation*}
d\le d_5(q):= \begin{cases}
d_4(q) & \text{if $q=p$,}\\
d_4(p,q) & \text{if $q=p^v$, $v\ge 2$,}
  \end{cases}
\qquad\text{or}\qquad d=\lfloor (q+2)/2\rfloor\, .
  \end{equation*}
   \end{theorem}
   \begin{proof} Suppose that $d>d_5(q)$. By means of some
   computations, $d_5(p,q)>d_4(q)$ and hence Lemma \ref{lemma2.3} holds true.
With the same notation as in the proof of that lemma, we can then
use the following two facts: $\epsilon=\nu\le q$, and $p|\epsilon$.
Actually, we will improve the latter one.
  \begin{claim}\label{claim2.1}
$\epsilon$ is a power of $p.$
   \end{claim}
Indeed, by $p|\epsilon$ and the $p$-adic criterion \cite[Cor.
1.9]{sv}, a necessary and sufficient condition for $\epsilon$ not
to be a power of $p$ is that $p\in\{2,3\}$ and $\epsilon=6$. If
this occurs, one can argue as in the previous proof and obtain the
following result:
$$
(5q-2)d^2-(q^2+15q-31)d+5(q+1)^2\le 0\, .
$$
{}From this,
$$
d\le G(q):=\frac{q^2+15q-31+
\sqrt{q^4-70q^3+203q^2-550q+1201}}{2(5q-12)}\, ,
$$
which is a contradiction as $G(q)<d_5(q)$.
   \begin{claim}\label{claim2.2}
\quad $\epsilon=q.$
   \end{claim}
The claim is certainly true for $q=p$. So, $q=p^v$, with $v\ge 2$.
If $\epsilon<q$, by Claim \ref{claim2.1} we have
$\epsilon\le q/p$. Thus, this fact together with
   $$
\deg(S_2)=(6+\nu)(2g-2)+(q^2+5)2d\ge 5\#\cX(\fq)=5(q+1)^2+5q(2g-2)\, ,
   $$
would yield
  $$
(5q-q/p-6)d^2-(2q^2+15q-3q/p-8)d+5(q+1)^2\le 0\, ,
  $$
and hence $d\le d_4(p,q)$, a contradiction.

Now from Claim \ref{claim2.2} and \cite[Cor. 2.6]{sv}, we have
   $$
q=\epsilon=\nu\le j_5(P_0)-1\le 2d-1\, ,
   $$
and Theorem \ref{thm2.1} follows from Corollary \ref{cor2.1}.
   \end{proof}

        \section{Maximal Hurwitz's curves}\label{s3}

In this section we give a necessary and sufficient condition for
$q$ in order that the Hurwitz curve $\cX_n$ defined by Eq.
(\ref{hurwitz}) be $\fq$-maximal.

   \begin{theorem}\label{thm3.1} The curve $\cX_n$ is $\fq$-maximal if and
only if (\ref{eq1.2}) holds.
   \end{theorem}
We first prove two lemmas.
  \begin{lemma}\label{lemma3.1} {\rm(\cite[p. 210]{carbonne-henocq})}
The Hurwitz curve $\cX_n$ is $\fp$-covered by the Fermat curve
$\cF_{n^2-n+1}$
   $$
U^{n^2-n+1}+V^{n^2-n+1}+W^{n^2-n+1}=0\, .
   $$
  \end{lemma}
  \begin{proof} Let $u=U/W$ and $v:=V/W$. Then the image of the morphism
$(u:v:1)\to (x:y:1)=(u^nv^{-1}:uv^{n-1}:1)$ is the curve defined by
$x^ny+y^n+x=0$. This proves the lemma.
   \end{proof}
  \begin{corollary}\label{cor3.1} Suppose that (\ref{eq1.2}) holds. Then
both curves $\cX_n$ and $\cF_{n^2-n+1}$ are $\fq$-covered by the
Hermitian curve of equation (\ref{hermitian}). In particular, both
are $\fq$-maximal.
  \end{corollary}
  \begin{proof} If (\ref{eq1.2}) holds, it is clear that $\cF_{n^2-n+1}$
is $\fq$-covered by the Hermitian curve. This property extends to
$\cX_n$ via the previous lemma. For both curves, the
$\fq$-maximality now follows from \cite[Prop. 6]{lachaud}.
  \end{proof}
  \begin{lemma}\label{lemma3.2} {\rm (\cite[p. 5249]{car-he-ri})} The
Weierstrass semigroup of $\cX_n$ at
the point $(0:1:0)$ is generated by the set
$S:= \{s(n-1)+1:s=1,\ldots,n\}$.
  \end{lemma}
  \begin{proof} Let $P_0:=(1:0:0)$, $P_1=(0:1:0)$, and $P_2=(0:0:1)$. Then
$\div(x)=nP_2-(n-1)P_1-P_0$ and $\div(y)=(n-1)P_0+P_2-nP_1$ so that
$$
\div(x^{s-1}y)=((n(s-1)+1)P_2+(n-s)P_0-(s(n-1)+1)P_1\, .
$$
This shows that $S$ is contained in the Weierstrass semigroup $H(P_1)$ at
$P_1$. In particular, $H(P_1)\supseteq \langle S\rangle$. Since
$\#(\N_0\setminus \langle S\rangle)=n(n-1)/2$ (see \cite{grant}), the
result follows.
  \end{proof}
{\em Proof of Theorem \ref{thm3.1}.} If (\ref{eq1.2}) holds, then $\cX_n$
is $\fq$-maximal by Corollary \ref{cor3.1}. Conversely, assume that
$\cX_n$ is $\fq$-maximal. Then $(q+1)P_1\sim
(q+1)P_2$ \cite[Lemma 1]{r-sti}, and the case $s=n$ in the proof of
Lemma \ref{lemma3.2} gives $(n^2-n+1)P_1\sim (n^2-n+1)P_2$. Therefore $d:={\rm
gcd}(n^2-n+1,q+1)$ belongs to $H(P_1)$. According to Lemma \ref{lemma3.2}
we have that $d=A(n-1)+B$ with $A\ge B\ge 1$. Now, there exists $C\ge 1$
such that $(A(n-1)+B)C=n^2-n+1$ and so $BC=D(n-1)+1$ for some $D\ge 0$.
Therefore, $AD(n-1)+A+BD=Bn$. We claim that $D=0$, otherwise the left side
of the last equality would be bigger than $Bn$. Then $B=C=1$ and so $A=n$;
i.e., $d=n^2-n+1$ and the proof is complete.

   \begin{corollary}\label{cor3.2} The curve $\cF_{n^2-n+1}$ in Lemma
\ref{lemma3.1} is $\fq$-maximal if and only if (\ref{eq1.2}) holds.
   \end{corollary}
   \begin{proof} If (\ref{eq1.2}) is satisfied, the result follows from
Corollary \ref{cor3.1}. Now if $\cF_{n^2-n+1}$ is $\fq$-maximal,
then $\cX_n$ is also $\fq$-maximal by Lemma \ref{lemma3.1} and
\cite[Prop. 6]{lachaud}. Then the corollary follows from Theorem
\ref{thm3.1}.
   \end{proof}
   \begin{remark}\label{rem3.0} For a given positive integer $n$,
we are led to look for a power $q$ of a prime $p$ such that
$q+1\equiv 0\pmod{m}$ with $m=n^2-n+1$. Since $m\not\equiv
0\pmod{p}$, and $p\not\equiv 0\pmod{m}$, a necessary and sufficient
condition for $q$ to have the requested property (\ref{eq1.2}) is
$p\equiv x\pmod{m}$, where $x$ is a 
solution of the congruence $X^w+1\equiv 0\pmod{m}$, and $w$ is
defined by $q=p^{\phi(m)v+w}$, $w\in\{1,2,\ldots,\phi(m)-1\}$; here 
$\phi$ denotes the Euler function.
   \end{remark}
Regarding specific examples, we notice that Carbonne and Henocq
\cite[Lemmes 3.3, 3.6]{carbonne-henocq} pointed out that $\cX_n$ is
$\fq$-maximal in the following cases:
  \begin{enumerate}
\item[\rm(1)] $n=3$, $q=p^{6v+3}$ and $p\equiv 3,5\pmod{7}$;
\item[\rm(2)] $n=4$, $q=p^{12v+6}$ and $p\equiv 2,6,7,11\pmod{13}$.
  \end{enumerate}
By using Theorem \ref{thm3.1} and Remark \ref{rem3.0} we have the
following result.
   \begin{corollary}\label{cor3.3}\begin{enumerate}
\item[\rm(1)] The curve $\cX_2$ is $\fq$-maximal if and only if
$q=p^{2v+1}$ and $p\equiv 2\pmod{3};$
\item[\rm(2)] The curve $\cX_3$ is $\fq$-maximal if and only if either
$q=p^{6v+1}$ and $p\equiv 6\pmod{7}$, or $q=p^{6v+3}$ and $p\equiv
3,5,6\pmod{7}$, or $q=p^{6v+5}$ and $p\equiv 6\pmod{7};$
\item[\rm(3)] The curve $\cX_4$ is $\fq$-maximal if and only if either
$q=p^{12v+1}$ and $p\equiv 12\pmod{13}$, or $q=p^{12v+2}$ and $p\equiv
5,8$, or $q=p^{12v+3}$ and $p\equiv 4,10, 12\pmod{13}$,
or $q=p^{12v+5}$ and $p\equiv 12\pmod{13}$, or
$q=p^{12v+6}$ and $p\equiv 2,5,6,7,8,11\pmod{13}$, or
$q=p^{12v+7}$ and $p\equiv 12\pmod{13}$, or
$q=p^{12v+9}$ and $p\equiv 4,10,12\pmod{13}$, or
$q=p^{12v+11}$ and $p\equiv 12\pmod{13}.$
   \end{enumerate}
   \end{corollary}
   \begin{corollary}\label{cor3.4} Let $n$ be a positive integer,
$m:=n^2-n+1$ and $p$ a prime.
\begin{enumerate}
    \item[\rm(1)] If $n=p^e$ with $e\ge 1$, then the curve $\cX_n$ is
$\fq$-maximal with $q=p^{\phi(m)v+3e}$.
    \item[\rm(2)] Let $p\equiv 3\pmod{4}$ and $n\equiv 0, 1\pmod{p}$ such
that $m$ is prime and that $m\equiv 3\pmod{4}$. Then $\cX_n$ is
$\fq$-maximal with $q=p^{(m-1)v+(m-1)/2}$.
\end{enumerate}
    \end{corollary}
    \begin{proof} Part (1) follows from the identity
$p^{3e}+1=(p^e+1)(p^{2e}-p^e+1)$ and Theorem \ref{thm3.1}.

To show (2), it is enough to check that $p^{(m-1)/2}+1\equiv
0\pmod{m}$. Recall that the Legendre symbol $(a/p)$ is defined by:
    \begin{equation*}
 (a/p)= \begin{cases}
   1   & \text{if $x^2\equiv a\pmod{p}$ has two solutions in $\Z_p$\, ,}
\\
  -1   & \text{if $x^2\equiv a\pmod{p}$ has no solution in $\Z_p$\, .}
 \end{cases}
    \end{equation*}
In our case, since $m\equiv 1\pmod{p}$, $(m/p)=1$. By the quadratic
reciprocity low
   $$
(m/p)(p/m)=(-1)^{((m-1)/2)((p-1)/2)}\, ,
   $$
from $(m/p)=1$ and $m\equiv 3\pmod{4}$ we get
$(p/m)=(-1)^{(p-1)/2}$. Now, as $p\equiv 3\pmod{4}$, we have that
$(p/m)=-1$. In other words, $p$ viewed as an element in
$\mathbf F_m$ is a non-square in $\mathbf F_m$. Since $-1$ is as
well a non-square in $\mathbf F_m$, it follows then that $p\equiv
(-1)u^2\pmod{m}$ with $u\in \Z$ such that $u\not\equiv 0\pmod{m}$.
Hence $p^{(m-1)/2}\equiv{(-1)}\pmod{m}$ as, in particular, $m$ is
odd and as $u^{m-1}\equiv 1\pmod m$.
   \end{proof}
   \begin{remark}\label{rem3.00} The hypothesis $m\equiv 3\pmod{4}$ in the
above corollary cannot be relaxed. In fact, for $n=4$ we have
$m=13$ but, according with Corollary \ref{cor3.3}, $\cX_4$ is no
$\mathbf F_{3^6}$-maximal.
   \end{remark}
   \begin{remark}\label{rem3.01} Let us assume the hypothesis in Corollary
\ref{cor3.4}(2) with $m$ not necessarily prime. In this case, to
study the congruence in (\ref{eq1.2}) we have to consider the
multiplicative group $\Phi_m$ of the units in $\Z_m$. This group
has order $\phi(m)$, and $p\in\Phi_m$ since $m\equiv 1\pmod{p}$.
Now suppose that $p$, as an element of $\Phi_m$, has even order
$2i$. Then $p^{2i}\equiv 1\pmod{m}$ and hence $(p^i+1)(p^i-1)\equiv
0\pmod{m}$. Since $p$ has order greater than $i$, we have that
$p^i-1\not\equiv0\pmod{m}$ unless both $p^i+1$ and $p^i-1$ are zero
divisors in $\Z_m$. If we assume that this does not happen, then
equivalence (\ref{eq1.2}) follows for $q=p^{\phi(m)v+i}$.
%
%This shows that we have to find out whether $p$ has even order in
%$\Phi_m$. It might be useful to take into account that $\Phi_m$ is direct
%product of cyclic groups. More precisely, if $m=p^{\alpha_1}_1\dots
%p^{\alpha_r}_r$ with primes $p_j$, then $\Phi_m$ is the product of cyclic
%groups of orders $p^{\alpha_j-1}_j(p_j-1)$, with $j=1,\dots,r$. To be
%continue.
   \end{remark}
   \begin{remark}\label{rem3.1} Let $p$ be a prime, $n:=p^eu$ with $e\ge
1$ and ${\rm gcd}(p,u)=1$. Assume $e\ge 2$ if $p=2$. Then
the Hurwitz curve $\cX_n$ as well as the curve $\cF_{n^2-n+1}$ are
non-classical with
respect to $\Sigma_1$. It is easy to see that $0,1$ and $p^e$ are their
$\Sigma_1$-orders.
    \end{remark}
%
%
%begin{remark}
%Is it true that the sequence $n^2-n+1$ with $n=1,\ldots,$ contains
%infinite
%primes
%$\equiv{3}\pmod{4}$? Does it hold true for ${n}\equiv{0}\pmod{p}$? We
%think
%that the
%answer is affermative to both questions. We are looking for results from
%number theory
%(Dirichlet theorem).
%\end{remark}
%
%

\section{On the maximality of generalized Hurwitz curves}\label{s4}

In this section we investigate the $\fq$-maximality of the non-singular
model of the so-called generalized Hurwitz curve $\cX_{n,\ell}$ of
equation
   $$
X^nY^\ell+Y^nZ^\ell+Z^nX^\ell=0\, ,
   $$
where $n\ge \ell\ge 2$ and $p={\rm char}(\fq)$ does not divide
$Q(n,\ell):= n^2-n\ell+\ell^2$. The singular points of $\cX_{n,\ell}$ are
$P_0:=(1:0:0)$, $P_1=(0:1:0)$, and $P_2=(0:0:1)$; each of them is
unibranched with $\delta$-invariant equal to $(n\ell-n-\ell+{\rm
gcd}(n,\ell))/2$. Therefore its genus $g$ (cf. \cite[Sec.
4]{bennama-carbonne} and \cite[Example 4.5]{beelen-pellikaan}) is equal to
$$
g=\frac{n^2-n\ell+\ell^2+2-3{\rm gcd}(n,\ell)}{2}\, .
$$
First we generalize Lemma \ref{lemma3.1}.
   \begin{lemma}\label{lemma4.1} The curve $\cX_{n,\ell}$ is $\fq$-covered
by the Fermat curve $\cF_{n^2-n\ell+\ell^2}$
   $$
U^{n^2-n\ell+\ell^2}+V^{n^2-n\ell+\ell^2}+W^{n^2-n\ell+\ell^2}=0\, .
   $$
   \end{lemma}
   \begin{proof} The curve $\cX_{n,\ell}$ is $\fq$-covered by
$\cF_{n^2-n\ell+\ell^2}$ via the morphism $(u:v:1)\to
(x:y:1):=(u^nv^{-m}:u^mv^{n-m}:1)$, where $u:=U/W$ and $v:=V/W$.
   \end{proof}
{}From this lemma and \cite[Prop. 6]{lachaud} we have the following.
   \begin{corollary}\label{cor4.1} The curve $\cF_{n^2-n\ell+\ell^2}$ in
the above lemma and the $\fq$-non-singular model of $\cX_{n,\ell}$ are
$\fq$-maximal provided that
\begin{equation}\label{eq4.1}
n^2-n\ell+\ell^2\equiv 0\pmod{(q+1)}\, .
   \end{equation}
   \end{corollary}
Now, we generalize Lemma \ref{lemma3.2} for any two coprime $n$ and
$\ell$. For $0 \leq i \leq 2$, let $Q_i$ be the unique point in the
non-singular model of $\cX_{n,\ell}$ lying over $P_i$.

   \begin{lemma}\label{lemma4.2} Suppose that ${\rm gcd}(n,\ell)=1$. Then
the Weierstrass semigroup $H(Q_1)$ at $Q_1$ is given by
\begin{equation}\label{eq4.2}
\{(n-\ell)s+nt: s,t\in \Z; t\ge 0\, -\frac{\ell}{n}t\le s\le
\frac{n-\ell}{\ell}t\}\, .
\end{equation}
   \end{lemma}
   \begin{proof} Let $x:=X/Z,y:=Y/Z$. It is not difficult to
see that $\div(x)=nQ_2-(n-\ell)Q_1-\ell Q_0$ and
$\div(y)=(n-\ell)Q_0+\ell Q_2-nQ_1$. Hence, for $s,t\in \Z$,
$$
\div(x^sy^t)=(ns+\ell t)Q_2+(-\ell s+(n-\ell)t)Q_0-((n-\ell)s+nt)Q_1\, ,
$$
and hence $(n-\ell)s+nt\in H(Q_1)$ provided that $ns+\ell t\ge 0$
and $-\ell s+(n-\ell)t\ge 0$. Let $H$ denote the set introduced in
(\ref{eq4.2}). Then $H\subseteq H(Q_1)$, and it is easily checked
that $H$ is a semigroup. By means of some computations we see that
$\#(\N\setminus H)= (n^2-n\ell+\ell^2-1)/2$, whence $H=H(Q_1)$
follows.
   \end{proof}
  \begin{remark}\label{rem4.1} The above Weierstrass semigroup $H(Q_1)$
was computed for $\ell=n-1$, and $(n,\ell)=(5,2)$ in
\cite{bennama-carbonne}.
  \end{remark}
We are able to generalize Theorem \ref{thm3.1} for certain
curves $\cX_{n,\ell}$.
   \begin{theorem}\label{thm4.1} Assume that ${\rm gcd}(n,\ell)=1$ and
that $Q:=Q(n,\ell)=n^2-n\ell+\ell^2$ is prime. Then $\cX_{n,\ell}$ is
$\fq$-maximal if and only if $(\ref{eq4.1})$ holds.
   \end{theorem}
  \begin{proof} The ``if" part follows from Corollary \ref{cor4.1} and
here we do not use the hypothesis that $Q$ is prime. For
the ``only if" part, we first notice that each $Q_i$ is $\fq$-rational.
Now the case $s=n-m$ and $t=m$ in
the proof of Lemma \ref{lemma4.2} gives $Q Q_2\sim Q Q_1$.
Therefore $d={\rm gcd}(Q,q+1)\in H(Q_1)$ because $(q+1)Q_1\sim
(q+1)Q_2$ \cite[lemma1]{r-sti}. As $1\not\in H(Q_1)$ and $Q$ is
prime, the result follows.
   \end{proof}
  \begin{corollary}\label{cor4.2} Let $n$, $\ell$ and $Q$ be as in Theorem
\ref{thm4.1}. Then the curve $\cF_{n^2-n\ell+\ell^2}$ in Lemma
\ref{lemma4.1} is $\fq$-maximal if
and only if (\ref{eq4.1}) holds.
   \end{corollary}
   \begin{proof} Similar to the proof of Corollary \ref{cor3.2}.
   \end{proof}
   \begin{remark}\label{rem4.2} There are infinitely many $n,\ell$ with
$n>\ell\ge 1$ such that $Q(n,\ell)$ is prime. In fact, for a prime $p'$
such that $p'\equiv 1\pmod{6}$, there exists such $n$ and $\ell$ so that
$p'=Q(n,\ell)$; see \cite[Remarque 4]{bennama-carbonne}.
   \end{remark}

\end{document}